\documentclass{lms}
\usepackage{amssymb}
\usepackage{amsmath,amsfonts}

\usepackage[ruled,vlined,linesnumbered]{algorithm2e}

\newtheorem{theorem}{Theorem}%[section]
\newtheorem{lemma}[theorem]{Lemma}
\newtheorem{proposition}[theorem]{Proposition}
\newnumbered{definition}{Definition}
\newnumbered{remark}{Remark}
\newnumbered{example}{Example}

\newcommand{\Q}{\mathbb{Q}}
\newcommand{\Z}{\mathbb{Z}}

\newcommand{\sfM}{{\sf M}}

\newcommand{\subgrp}[1]{\ensuremath{\langle{#1}\rangle}}

\newcommand{\QQ}{\Q}
\newcommand{\ZZ}{\Z}

\newcommand{\FF}{\mathbb{F}}
\newcommand{\FFbar}{\overline{\mathbb{F}}}
\newcommand{\EC}{\mathcal{E}}
\newcommand{\ECK}{\widetilde{\mathcal{E}}}
\newcommand{\phiK}{\widetilde{\phi}}
\newcommand{\conj}[2][\sigma]{{{}^{#1}{#2}}}
\newcommand{\pconj}[1]{{\conj{#1}}}
\newcommand{\dualof}[1]{{{#1}^{\dagger}}}
\newcommand{\End}{\mathrm{End}}
\newcommand{\Aut}{\mathrm{Aut}}
\newcommand{\isogeny}{\ensuremath{\vartheta}}

\title[Computing cardinalities of \(\Q\)-curve reductions]{Computing cardinalities of $\Q$-curve reductions over finite fields}
\author{Fran\c{c}ois Morain, Charlotte Scribot, and Benjamin Smith}

\begin{document}

\maketitle

\begin{abstract}
    We present a specialized point-counting algorithm
    for a class of elliptic curves over \(\FF_{p^2}\)
    that includes reductions of quadratic \(\QQ\)-curves modulo inert primes
    and, more generally,
    any elliptic curve over \(\FF_{p^2}\)
    with a low-degree isogeny to its Galois conjugate curve.
    These curves have interesting cryptographic applications.
    Our algorithm is a variant of the Schoof--Elkies--Atkin (SEA)
    algorithm, but with a new, lower-degree endomorphism in place of Frobenius.
    While it has the same asymptotic asymptotic complexity as SEA,
    our algorithm is much faster in practice.
\end{abstract}

%%%%% S
\section{Introduction}

Computing the cardinalities of the groups of rational points
on elliptic curves over finite fields is a fundamental algorithmic
challenge in computational number theory, 
and an essential tool in elliptic curve cryptography.
Over finite fields of large characteristic,
the best known algorithm is the Schoof--Elkies--Atkin (SEA)
algorithm~\cite{Schoof95}. A lot of work has been put into optimizing the
computations for prime fields of large characteristic (see
\cite{Sutherland13} for the most recent record). Many of these
improvements carry over to the case of more general finite fields.
In this article we define a specialized, faster SEA algorithm
for a class of elliptic curves over \(\FF_{p^2}\)
that have useful cryptographic applications.
These curves have low-degree inseparable endomorphisms
that can be used to accelerate scalar multiplication
in elliptic curve cryptosystems~\cite{Smith13,Smith15};
here, we use those endomorphisms to accelerate point counting.
Going beyond cryptography,
this class of curves also includes reductions of quadratic \(\QQ\)-curves
modulo inert primes, so our algorithm may be useful for studying
these curves.

Let \(q\) be a power of a prime \(p > 3\)
(in our applications, \(q = p^2\) and \(p\) is large).
Let
\[
    \pconj{(\cdot)} : x \longmapsto x^p
\]
be the Frobenius automorphism on \(\FF_q\).
We extend the action of Frobenius to polynomials over \(\FF_q\)
by acting on coefficients,
and thus to curves over \(\FF_q\) by acting on their defining equations:
for example, an elliptic curve \(\EC/\FF_q\) and its Galois conjugate curve
\(\pconj{\EC}/\FF_q\) would be defined by
\[
    \EC: y^2 = x^3 + Ax + B
    \quad
    \text{and}
    \quad
    \pconj{\EC}: y^2 = x^3 + A^px + B^p
    \ .
\]

If \(\EC/\FF_q\) is an elliptic curve,
then there is a \(p\)-isogeny
\(\pi_p: \EC\to\pconj{\EC}\)
defined by \(\pi_p : (x,y) \mapsto (x^p, y^p)\).
If \(q = p^n\),
then composing 
\(\pi_p\), \(\pconj{\pi_p}\), \dots,
\(\conj[\sigma^{n-1}]{\pi_p}\)
yields the Frobenius endomorphism 
\(\pi_q : (x,y)\mapsto (x^q,y^q)\) of \(\EC\).
Being an endomorphism, \(\pi_q\) has a characteristic polynomial
\[
    \label{eq:Frobenius-poly}
    \chi_{\pi_{q}}(T)
    = 
    T^2 - t_\EC T + q
\]
such that \(\chi_{\pi_q}(\pi_q) = [0]\) in \(\End(\EC)\);
the \emph{trace} \(t_\EC\) satisfies the Hasse bound 
\[
    \label{eq:Hasse--Weil}
    |t_\EC| \le 2\sqrt{q}
    \ .
\]
Knowing the cardinality of \(\EC(\FF_q)\)
is equivalent to knowing the trace \(t_\EC\),
since
\[
    \#\EC(\FF_{q}) 
    =
    \chi_{\pi_{q}}(1) 
    = 
    q + 1 - t_\EC 
    \ .
\]

Schoof's point counting algorithm~\cite{Schoof85} 
determines $t_\ell := t_\EC \pmod{\ell}$
for small primes $\ell\not= p$ 
by examining the action of $\pi_q$ 
on \(\EC[\ell]\), the $\ell$-torsion subscheme of $\EC$:
we have
\[
    \pi_q^2(P) - [t_\ell] \pi_q(P) + [q\bmod\ell] P = 0
    \quad
    \text{ for } P \in \EC[\ell] 
    \ .
\]
If we construct a general \(P\) as detailed in~\S\ref{sec:background},
then finding \(t_\ell\) boils down to a series of polynomial operations 
modulo the \(\ell\)-th division polynomial \(\Psi_\ell\).
Schoof's algorithm tests these relations until $\prod \ell > 4
\sqrt{q}$, and then deduces \(t_\EC\) from the \(t_\ell\)
using the Chinese Remainder Theorem (CRT). 
To completely determine \(t_\EC\) 
we need to compute \(t_\ell\) for $O(\log q)$ primes \(\ell\),
the largest of which is in \(O(\log q)\); 
fast polynomial evaluations add some more $O(\log q)$ factors, 
and the final cost is $O(\log^8 q)$ with classical arithmetic (or
$O(\log^6 q)$ with fast arithmetic).
This basic algorithm was subsequently improved by Atkin and Elkies;
the resulting SEA algorithm (see~\S\ref{sec:SEA})
is now the standard point-counting
algorithm for elliptic curves over large characteristic fields.

In this article,
we present an algorithm that was designed to compute \(\#\EC(\FF_{p^2})\) 
when \(\EC\) is the reduction
of a low-degree quadratic \(\QQ\)-curve modulo an inert prime.
In fact, our algorithm applies to a larger class of curves over finite fields,
which we will call \emph{admissible curves}.

First, 
recall that every \(d\)-isogeny \(\isogeny: \EC \to \EC'\) 
has a \emph{dual} \(d\)-isogeny \(\dualof{\isogeny}: \EC' \to \EC\)
such that \(\dualof{\isogeny}\isogeny = [d]_{\EC}\)
and \(\isogeny\dualof{\isogeny} = [d]_{\EC'}\).
Also, \(\sigma\) acts on isogenies
by \(p\)-th powering the coefficients of their defining polynomials;
so
every isogeny \(\isogeny : \EC \to \EC'\) 
has a Galois conjugate isogeny
\(\pconj{\isogeny}: \pconj{\EC} \to \pconj{\EC'}\).

\begin{definition}
    \label{def:admissible}
    Let \(d\) be a squarefree integer with \(p\nmid d\).
    An elliptic curve \(\EC/\FF_{p^2}\) is \emph{\(d\)-admissible}
    if it is equipped with a \(d\)-isogeny
    \[
        \label{eq:conjugate-duality}
        \phi: \EC \longrightarrow \pconj{\EC}
        \quad
        \text{ such that }
        \quad
        \pconj{\phi}
        =
        \epsilon\dualof{\phi} 
        \quad \text{where} \quad 
        \epsilon = \pm 1
        \ .
    \]
    Composing 
    \(\pi_p: \EC \to \pconj{\EC}\)
    with \(\pconj{\phi} : \pconj{\EC} \to \EC\),
    we obtain
    the \emph{associated endomorphism}
    \[
        \psi := \pconj{\phi}\circ\pi_p
        \in \End(\EC)
    \]
    of degree \(dp\).
    Note that the requirement \(p\nmid d\)
    implies that both \(\phi\) and \(\pconj{\phi}\)
    are separable.
\end{definition}

We are particularly interested in curves that are \(d\)-admissible 
for small values of \(d\).
% FIXME: John Cremona said that there are (infinitely?) many degrees of
% separable isogenies to the Galois conjugate.  Can we find a reference
% for this?
When \(d\) is extremely small
the associated endomorphism can be evaluated very efficiently,
and thus used to accelerate scalar multiplication on \(\EC\)
for more efficient implementations of elliptic curve cryptosystems
(as in~\cite{Smith13}, \cite{GuillevicIonica13}, \cite{CoHiSm14},
\cite{Smith15}, and~\cite{CostelloLonga15}).
Constructing cryptographically secure curves 
equipped with efficient endomorphisms
is one major motivation for our algorithm;
the other is the principle that the presence of special structures
demands the use of a specialized algorithm.

From a practical point of view, suitable modifications of the 
SEA algorithm gives us a very fast probabilistic solution to the 
point counting problem for admissible curves.
The essential idea is to use SEA
with the associated endomorphism \(\psi\)
in place of \(\pi_q\).
While the asymptotic complexity of our algorithm is the same as for the
unmodified SEA algorithm when \(d\) is fixed, there are some important 
improvements in the big-\(O\) constants.  
Asymptotically, when \(d\) is small, 
our algorithm runs four times faster than SEA
(and even faster for smaller \(p\)).

It is not hard to see that of the \(p^2\) isomorphism classes of
elliptic curves over \(\FF_{p^2}\), only \(O(p)\) classes
correspond to \(d\)-admissible curves for any \emph{fixed} \(d\).
But while \(d\)-admissible curves with small \(d\) may be relatively rare,
they appear naturally ``in the wild'' 
as reductions of quadratic \(\QQ\)-curves
of degree \(d\)
(elliptic curves over quadratic number fields that are \(d\)-isogenous to their
Galois conjugates) modulo inert primes.
For some small \(d\), these \(\QQ\)-curves occur in one-parameter families;
so our algorithm allows the reductions of these families modulo suitable
primes to be rapidly searched for cryptographic curves.
We explain this further in~\S\ref{sec:QQ-curves}.

%%%%% S
\section{%%%%%%%%%%%%%%%%%%%%%%%%%%%%%%%%%%%%%%%%%%%%%%%%%%%%%%%%%%%%%%%%%%%%%%%
    Computing with isogenies
}%%%%%%%%%%%%%%%%%%%%%%%%%%%%%%%%%%%%%%%%%%%%%%%%%%%%%%%%%%%%%%%%%%%%%%%%%%%%%%%
\label{sec:background}

We begin by recalling some standard results on isogenies,
fixing notation and basic complexities in the process. A classical
reference for all this is \cite{GaGe99}.

First, let ${\sfM}(n)$ denote the cost in \(\FF_q\)-operations
(multiplications) of multiplying two polynomials of degree~$n$.
Traditional multiplication gives
${\sfM}(n) = O(n^2)$;
fast multiplication gives $\widetilde{O}(n)$.
Dividing a degree-\(2n\) polynomial by 
a degree-\(n\) polynomial costs $O({\sfM}(n))$ \(\FF_q\)-operations;
the extended GCD of two degree-\(n\) polynomials
can be computed in $O({\sfM}(n) \log n)$ \(\FF_q\)-operations.
The number of roots in \(\FF_q\)
of a degree-\(n\) polynomial \(F\) over \(\FF_q\)
is equal to \(\deg\textsc{Gcd}(x^q-x,F(x))\),
which we can compute in \(O((\log q)\sfM(n))\) \(\FF_q\)-operations
if \(n \ll q\)
(this is dominated by the cost of computing \(x^q \bmod F\);
see Appendix~\ref{sec:technical-lemmas}).

We will make extensive use of \emph{modular composition}: 
if $F$, \(G\), and \(H\) are polynomials over \(\FF_q\)
with \(\deg F = n\), $\deg G < n$, and $\deg H < n$, 
then we can compute $(G \circ H) \bmod F$
in $O(n^{1/2} {\sfM}(n) + n^{(\omega+1)/2})$ \(\FF_q\)-operations,
where $2 \le \omega \le 3$ is the constant for linear algebra.
Using the method of~\cite{KaSh98},
the cost in \(\FF_q\)-operations of performing $r$ modular compositions 
with the same $H$ and \(F\) is 
\[
    \mathcal{C}_r(n) 
    := 
    O(r^{1/2} n^{1/2} {\sfM}(n) + r^{(\omega-1)/2} n^{(\omega+1)/2})
    \ .
\]
\iffalse
    \item \textbf{Iterated modular composition}: 
        given $k \in O(n)$, 
        we want to compute $h,h(h) \bmod F, \dots,h^{(k)} \bmod F,$
        assuming $F$ divides $F(h)$.  
        Using the algorithm of~\cite[Lemma 4]{KaSh98}: 
        if $$H_1=h,H_2=h(h) \bmod F, \dots,H_j=h^{(j)} \bmod F$$
        are known, then we deduce
        $$H_{j+1}=H_1(H_j) \bmod F, \dots,H_{2j}=H_j(H_j) \bmod F.$$
        We repeat this scheme for $j=1,\dots,2^{\lceil \log(k) \rceil}$.  
        The total cost is up a constant that of the last step, {\it i.e.},
        $O(\mathcal{C}_k(n))$.
\fi

We will always work with elliptic curves \(\EC/\FF_q\)
using their Weierstrass models,
\[
    \EC: y^2 = f_\EC(x) ,
    \text{ where } f_\EC \text{ is a monic cubic over } \FF_q\ .
\]
For \(m > 0\),
the \(m\)-th \emph{division polynomial}
\(\Psi_m(x)\) is the polynomial in \(\FF_q[x]\)
whose roots are precisely the
\(x\)-coordinates of the points in \(\EC[m](\FFbar_{q})\).

If \(\ell\) is a prime,
then the level-\(\ell\) modular polynomial 
\(\Phi_\ell(J_1,J_2)\) has degree \(\ell+1\) 
in both \(J_1\) and \(J_2\),
and is defined over \(\ZZ\).
If \(\Phi_\ell(j_1,j_2) = 0\) for some \(j_1\) and \(j_2\) in
\(\FF_{q}\),
then there is an \(\FF_{q}\)-rational \(\ell\)-isogeny
between the curves with \(j\)-invariants \(j_1\) and \(j_2\)
(possibly after a twist).
In particular,
if we fix an elliptic curve \(\EC/\FF_q\),
then the roots of \(\Phi_\ell(j(\EC),x)\) in \(\FF_q\)
correspond to (the isomorphism classes of) the curves that are
\(\ell\)-isogenous to \(\EC\) over \(\FF_q\).

We will need explicit forms for \(d\)-isogenies where \(d\) is
squarefree and prime to \(p\).
Every such isogeny can be expressed as a composition of
at most one \(2\)-isogeny with at most one odd-degree cyclic isogeny
over \(\FF_q\).
If \(\isogeny\) is a \(2\)-isogeny,
then it is defined by a rational map
\begin{equation}
    \label{eq:2-isogeny-def}
    \isogeny:
    (x,y) 
    \longmapsto 
    \left(
        \frac{N(x)}{D(x)},
        y\frac{M(x)}{D^2(x)}
    \right)
\end{equation}
where \(N\), \(M\), and \(D\)
are polynomials over \(\FF_q\)
with \(\deg N = \deg M = 2\) and \(D = x-x_0\) where $x_0$ is the abscissa
of a $2$-torsion point.
If \(\isogeny\) is a \(d\)-isogeny
where \(d\) is odd, squarefree, and prime to \(p\),
then \(\isogeny\)
is defined by a rational map 
\begin{equation}
    \label{eq:d-isogeny-def}
    \isogeny:
    (x,y)
    \longmapsto
    \left(
        \frac{N(x)}{D^2(x)},
        y\frac{M(x)}{D^3(x)}
    \right)
\end{equation}
where \(N\), \(M\), and \(D\)
are polynomials over \(\FF_q\)
with \(\deg N = \deg M = d\) and \(\deg D = (d-1)/2\).

In both cases, 
the polynomial \(D(x)\) cuts out the kernel of \(\isogeny\),
in the sense that \(D(x(P)) = 0\) if and only if \(P\)
is a nontrivial element of \(\ker\isogeny\);
we call \(D\) the \emph{kernel polynomial} of \(\isogeny\).
We suppose we have a subroutine
\(\textsc{KernelPolynomial}(\ell,\EC,j_1)\)
which, given \(\EC\) and \(j_1 = j(\EC_1)\)
such that there exists an \(\ell\)-isogeny
\(\isogeny:\EC\to\EC_1\) over \(\FF_q\),
computes the kernel polynomial \(D\) of \(\isogeny\)
and the isogenous curve \(\EC_1\)
in $O(\ell^2)$ \(\FF_q\)-operations
(using the fast algorithms in~\cite{BMSS08}).

The algorithms in this article
examine the actions of endomorphisms
on \(\ker\isogeny\),
where \(\isogeny\) is either \([\ell]\) or an \(\ell\)-isogeny,
for a series of small primes \(\ell\).
The key is to define a symbolic element of \(\ker\isogeny\).
First, we
compute the kernel polynomial \(D\) of \(\isogeny\)
(note that \(D = \Psi_\ell\) if \(\isogeny = [\ell]\));
then,
we can construct a symbolic point \(P\) of \(\ker\isogeny\) 
as
\[
    P := (X,Y) 
    \in 
    \EC\left(\FF_q[X,Y]/(Y^2 - f_\EC(X),D(X))\right)
    \ .
\]
We reduce the coordinates of points in \(\subgrp{P}\) 
modulo
\(D(X)\) and \(Y^2 - f_\EC(X)\) after each operation,
so elements of \(\subgrp{P}\)
have a canonical form \(Q = (Q_x(X),YQ_y(X))\)
with \(\deg Q_x,\deg Q_y < \deg D\).

Let \(e = \deg D\);
then we can compute \(Q_1 + Q_2\) 
for any \(Q_1\) and \(Q_2\) in \(\subgrp{P}\)
in $O({\sfM}(e) \log e)$ \(\FF_q\)-operations,
using the standard affine Weierstrass addition formul\ae{}.
We can therefore compute \([m]Q\)
for any \(m\) in \(\ZZ\) and \(Q\) in \(\subgrp{P}\)
in $O((\log m) {\sfM}(e) \log e)$ \(\FF_q\)-operations,
using a binary method.
We let \(\textsc{DiscreteLogarithm}(Q_1,Q_2)\)
be a subroutine which returns 
the discrete logarithm of \(Q_2\) to the base \(Q_1\), 
where both points are in \(\subgrp{P}\),
in $O(\sqrt{e} {\sfM}(e))$ \(\FF_q\)-operations
(using the approach in~\cite{GaMo06};
in some cases we can do better~\cite{MiMoSc07}).

\begin{lemma}
    \label{lemma:costs}
    Let \(P = (X,Y)\)
    in \(\EC(\FF_q[X,Y]/(Y^2 - f_\EC(X),D(X)))\),
    and let \(e = \deg D\).
    Then for any \(Q\) in \(\subgrp{P}\), we can
    \begin{enumerate}
        \item 
            compute $\pi_p(P) = (X^p, Y^p)$ 
            in $O((\log p) \sfM(e))$ \(\FF_q\)-operations;
        \item 
            compute \(\pi_p(Q)\),
            given $\pi_p(P)$,
            in $O((\log p) {\sfM}(e))$ \(\FF_q\)-operations;
        \item 
            compute \(\phi(Q)\),
            where \(\phi\) is a \(2\)-isogeny
            (as in~\eqref{eq:2-isogeny-def})
            in $O({\sfM}(e))$ \(\FF_q\)-operations;
        \item
            compute \(\phi(Q)\),
            where \(\phi\) is a \(d\)-isogeny
            with \(d\) odd, squarefree, and prime to \(p\)
            (as in~\eqref{eq:d-isogeny-def})
            in $O({\sfM}(d) + \mathcal{C}_3(e))$ \(\FF_q\)-operations.
    \end{enumerate}
\begin{proof}
    See Appendix~\ref{sec:technical-lemmas}.
\end{proof}
\end{lemma}

%%%%% S
\section{Atkin, Elkies, and volcanic primes}
\label{sec:prime-classes}

Given an elliptic curve \(\EC/\FF_q\),
we split the primes \(\ell \not= p\) into three classes: 
Elkies, Atkin, and volcanic.
The volcanic primes fall in two sub-classes:
floor-volcanic and upper-volcanic.
This classification reflects the structure of the \(\ell\)-isogeny
graph near \(\EC\), 
which reflects the factorization of \(\Phi_\ell(j(\EC),x)\).
The facts stated below without proof all follow immediately from well-known
observations of Atkin for general ordinary elliptic curves over \(\FF_q\) 
(cf.~\cite[Prop.~6.2]{Schoof95}).

Recall that the discriminant of \(\chi_{\pi_{q}}\) 
is \( \Delta_{\pi_q} := t_\EC^2 - 4q  < 0\).
We say that \(\ell\) is \textbf{volcanic}
if \(\ell\) divides \(\Delta_{\pi_q}\).
A volcanic prime \(\ell\) is 
\textbf{floor-volcanic} if 
\[
    \label{eq:volcanic-factorization-I}
    \Phi_\ell(x,j(\EC))
    =
    (x - j_1)h(x)
    \ ,
\]
where \(h\) is an \(\FF_{q}\)-irreducible polynomial of degree \(\ell\),
or \textbf{upper-volcanic}
if
\[
    \label{eq:volcanic-factorization-II}
    \Phi_\ell(x,j(\EC))
    =
    \prod_{i=1}^{\ell+1}(x - j_i)
\]
with each \(j_i\) in \(\FF_q\).
In each case,
the roots \(j_i\) 
are the \(j\)-invariants of the elliptic curves \(\EC_i\) 
that are \(\ell\)-isogenous to \(\EC\) over \(\FF_{q}\)
(up to isomorphism).

We say that \(\ell\) is \textbf{Elkies}
if \(\Delta_{\pi_q}\) is a nonzero square modulo \(\ell\).
Equivalently, 
\(\ell\) is Elkies
if 
\[
    \label{eq:Elkies-factorization}
    \Phi_\ell(x,j(\EC))
    =
    (x - j_1) (x - j_2) \prod_{i=1}^{(\ell-1)/e} h_i(x) 
    \ ,
\]
where \(j_1\) and \(j_2\) are in \(\FF_q\)
and the \(h_i\) are \(\FF_q\)-irreducible polynomials,
all of the same degree \(e > 1\),
with \(e\mid (\ell-1)\).
In this case, there exist \(\FF_q\)-rational \(\ell\)-isogenies
\(\isogeny_1: \EC\to\EC_1\) and \(\isogeny_2: \EC\to\EC_2\)
such that \(j(\EC_i) = j_i\),
and the \(\ell\)-torsion decomposes as
\(\EC[\ell] = \ker\isogeny_1\oplus\ker\isogeny_2\).

We say that \(\ell\) is \textbf{Atkin}
if \(\Delta_{\pi_q}\) is not a square modulo \(\ell\).
Equivalently,
\(\ell\) is Atkin if
\[
    \label{eq:Atkin-factorization}
    \Phi_\ell(x,j(\EC))
    =
    \prod_{i=1}^{(\ell+1)/e}
    h_i(x)
    \ ,
\]
where the \(h_i\) are all 
irreducible polynomials of the same degree \(e > 1\),
with \(e \mid (\ell+1)\).
Since \(\Phi_\ell(x,j(\EC))\) has no roots in \(\FF_{q}\),
there are no elliptic curves 
\(\ell\)-isogenous to \(\EC\)
over \(\FF_{q}\).

We can determine the class of a prime \(\ell\)
by finding out how many roots \(\Phi_\ell(j(\EC),x)\) has in~\(\FF_q\).
We define a subroutine
\(\textsc{EvaluatedModularPolynomial}(\ell,\EC)\),
which computes \(\Phi_\ell(j(\EC),x)\)
in $O(\ell^3 (\log \ell)^3 \log\log \ell)$ bit operations (under
the GRH)
using the method of~\cite{Sutherland13}, assuming $\log q = \Theta(\ell)$.
(Note that in practice, one generally uses precomputed modular
polynomials over~$\Z$.)

The number of roots is 
the degree of 
\(J = \textsc{Gcd}(x^q-x,\Phi_\ell(j(\EC),x))\),
which we compute
at a further cost of \(O((\log q)\sfM(\ell))\)
\(\FF_q\)-operations.
We may then want one of these roots, if any exist;
we therefore define a subroutine
\(\textsc{OneRoot}(J)\)
which finds a single root of \(J\).
At worst, in the upper-volcanic case,
this requires 
\(O((\log q)M(\deg J)\log\deg J) = O((\log q)M(\ell)\log\ell)\) 
\(\FF_q\)-operations;
at best, in the lower-volcanic and Elkies cases
(where \(J\) is linear and quadratic, respectively),
\(\textsc{OneRoot}(J)\) costs \(O(1)\) \(\FF_q\)-operations.

%%%%% S
\section{The SEA algorithm}
\label{sec:SEA}

Algorithm~\ref{alg:SEA}
presents a basic version of the SEA algorithm.
The main loop computes \(t_\ell := t_\EC \pmod{\ell}\)
for a series of small primes \(\ell\);
then we recover \(t_\EC\) from the \(t_\ell\)
via the CRT.

The complexity of Algorithm~\ref{alg:SEA} (and
Algorithm~\ref{alg:AdmissibleTrace} below)
depends on the number of non-Atkin primes less than a given bound.
The standard (and na\"ive) heuristic on prime classes is
to suppose that
the number of Atkin and non-Atkin primes \(\ell\) less than \(B\) 
for a given \(\EC/\FF_q\)
is approximately equal when \(B \sim \log q\), as \(q \to \infty\).
In particular, this means that 
\(O(\log q)\) non-Atkin \(\ell\) suffice to determine \(t_\EC\), 
and the largest such \(\ell\) is in \(O(\log q)\).
While the standard heuristic holds on the average,
it is known to fail for some \(\EC\); 
Galbraith and Satoh have shown (under the GRH) that for some \(\EC/\FF_p\)
we may need to use non-Atkin \(\ell\) as large as \(O(\log^{2+\epsilon}p)\)
(see~\cite[App.~A]{Satoh02}).
We refer the reader to~\cite{ShSu14} and~\cite{ShSu15} 
for further details and discussion.

\begin{proposition}
    If \(\EC/\FF_q\) is an elliptic curve,
    then
    under the standard heuristic on prime classes,
    Algorithm~\ref{alg:SEA}
    computes \(t_\EC\) in 
    \(\widetilde{O}(\log^3 q)\) expected \(\FF_q\)-operations
    (that is \(\widetilde{O}(\log^4q)\) expected bit operations,
    using fast arithmetic).
%% L = O(log q); dominant term is computation of X^q mod Phi, so
%% total time is (log q) (log q) M(log q) operations over FF_q.
\begin{proof}
    The main loop computes a set \(\mathcal{T}\) 
    of pairs \((t_\ell := t_\EC \pmod{\ell},\ell)\)
    with \(\prod_\mathcal{T} \ell > 4\sqrt{q}\).
    We then recover \(t_\EC\) from \(\mathcal{T}\) 
    via an explicit CRT. % (Line~21).
    Our procedure for computing \(t_\ell\)
    depends on the class of \(\ell\),
    which we determine using the method at the end of~\S\ref{sec:prime-classes}
    (Lines~\ref{line:SEA:class}, \ref{line:SEA:is-volcanic}, 
    and~\ref{line:SEA:is-Elkies}).
        
    If \(\ell\) is volcanic
    (Lines~\ref{line:SEA:volcanic-start} to~\ref{line:SEA:volcanic-end}), 
    then \(\ell\mid\Delta_{\pi_q}\),
    so \(t_\ell = 0\) or \(t_\ell \equiv \pm2\sqrt{q}\pmod{\ell}\).
    %according to~\eqref{eq:Atkin-trace-relation}.
    We distinguish between the three cases by
    comparing \(\pi_q(P)\) with \(\pm[\sqrt{q}\bmod\ell]P\)
    for a generic element \(P\)
    of the kernel 
    of the rational \(\ell\)-isogeny 
    corresponding to one of the roots of \(\Phi_\ell(j(\EC),x)\).
    
    If \(\ell\) is Elkies
    (Lines~\ref{line:SEA:Elkies-start} to~\ref{line:SEA:Elkies-end}), 
    then \(\EC[\ell]\)
    decomposes as a direct sum
    \((\ker\isogeny_1)\oplus(\ker\isogeny_2)\)
    of \(\ell\)-isogeny kernels;
    \(\pi_q\)
    acts as multiplication by eigenvalues \(\lambda_1\) and \(\lambda_2\) 
    on \(\ker\isogeny_1\) and \(\ker\isogeny_2\), respectively,
    with \(\lambda_1\lambda_2 \equiv q \pmod{\ell}\),
    so \(t_\ell \equiv \lambda_1 + q/\lambda_1 \pmod{\ell}\);
    and we can determine \(\lambda_1\)
    by solving the discrete logarithm problem
    \(\pi_q(P) = [\lambda_1]P\)
    for a symbolic point \(P\)
    of \(\ker\isogeny_1\).

    If \(\ell\) is Atkin,
    then we skip it completely
    and do not compute \(t_\ell\)
    (see the discussion in~\S\ref{sec:complements}).

    In terms of \(\FF_q\)-operations,
    determining the class of \(\ell\)
    costs 
    \(O(\ell^2(\log\ell)^3\log\log\ell + \log q \sfM(\ell))\);
    computing \(t_\ell\) then costs
    \(O(
        (\log q + \log\ell)\sfM(\ell)\log\ell
        + 
        \ell^{(\omega+1)/2}
    )\) 
    for volcanic \(\ell\), 
    and
    \(O((\log p + \ell^{1/2})\sfM(\ell) + \ell^{(\omega+1)/2})\)
    for Elkies \(\ell\).
    % cost: O( #primes ) * O(max class + max Elkies + max volcanic)
    % = O(log q)*(O~(log^3 q) + O(log q M(log q)) + O(log qM(log q)))
    The standard heuristic on prime classes
    tells us that we will try \(O(\log q)\) primes \(\ell\),
    and that the largest \(\ell\) are in \(O(\log q)\);
    so the total cost of the algorithm is
    \(\widetilde{O}(\log^3 q)\), as claimed.
\end{proof}
\end{proposition}

\begin{algorithm}
    \caption{\textsc{SEATrace}}
    \label{alg:SEA}
    \KwIn{An elliptic curve \(\EC/\FF_q\), where \(q = p^n\) with \(p\)
    large}
    \KwOut{The trace of Frobenius of \(\EC\)}
    % Choose a set \(\mathcal{L}\) of primes \(\ell\not=p\)
    % such that \(\prod_{\ell\in\mathcal{L}}\ell \ge 4p\) 
    % \;
    \(\mathcal{T} \gets \{ \}\)
    \tcp*{\(\mathcal{T}\) will contain the pairs \((t_\EC\pmod{\ell},\ell)\)}
    \(M \gets 1\)
    \tcp*{After each iteration, \(t_\EC\) is known modulo \(M\)}
    \(\ell \gets 1\)
    \;
    \While{\(M \le 4\sqrt{q}\)}{
        \(\ell \gets \textsc{NextPrime}(\ell)\)
        \;
        \(J \gets \textsc{Gcd}(x^q-x,\textsc{EvaluatedModularPolynomial}(\ell,\EC))\)
        \; \label{line:SEA:class}
        \If(\tcp*[f]{\(\ell\) is volcanic}){\(\deg J = 1\) or \(\ell+ 1\)}{
            \label{line:SEA:is-volcanic}
            \If(\tcp*[f]{\(s = p\) for \(q = p^2\)}){\(q\) has a square
            root \(s\) modulo \(\ell\)}{
                \(F \gets
            	\textsc{KernelPolynomial}(\ell,\EC,\textsc{OneRoot}(J))\)
            	\; \label{line:SEA:volcanic-start}
            	\(P \gets (X,Y)\) 
            	in \(\EC(\FF_q[X,Y]/(Y^2 - f_\EC(X),F(X)))\)
            	\;
                \(Q_1 \gets \pi_p(P)\)
                ;
                \(Q_2 \gets \pi_p(Q_1)\)
                ;
                \(Q_3 \gets [s]P\)
                \;
                \(
                    t_\ell 
                    \gets
                    \begin{cases}
                        -2s & \text{if } Q_2 = Q_3 \ ; \\
                        2s  & \text{if } Q_2 = -Q_3 \\
                        0   & \text{otherwise} \ ;
                    \end{cases}
                \)
                \;
            }
            \lElse{\(t_\ell \gets 0\)}
            \(\mathcal{T} \gets \mathcal{T} \cup \{(t_\ell,\ell)\}\)
            ;
            \(M \gets \ell M\)
            \; \label{line:SEA:volcanic-end}
        }
        \ElseIf(\tcp*[f]{\(\ell\) is Elkies}){\(\deg J = 2\)}{
            \label{line:SEA:is-Elkies}
            \(F \gets
            \textsc{KernelPolynomial}(\ell,\EC,\textsc{OneRoot}(J))\)
            \; \label{line:SEA:Elkies-start}
            \(P \gets (X,Y)\) 
            in \(\EC(\FF_q[X,Y]/(Y^2 - f_\EC(X),F(X)))\)
            \;
            \(Q_1 \gets \pi_p(P)\)
            ;
            \(Q_2 \gets \pi_p(Q_1)\)
            \;
            \(t_\ell \gets \lambda + q/\lambda \pmod{\ell}\)
            \textbf{where}
            \(\lambda = \textsc{DiscreteLogarithm}(P,Q_2)\)
            \;
            \(\mathcal{T}\gets\mathcal{T}\cup\{(t_\ell,\ell)\}\)
            ;
            \(M \gets \ell M\)
            \; \label{line:SEA:Elkies-end}
        }
%        \Else(\tcp*[h]{\(\ell\) is Atkin}){
%            \(F \gets \textsc{DivisionPolynomial}(\ell,\EC)\)
%            \;
%            \(P \gets (X,Y)\) 
%            in \(\EC(\FF_q[X,Y]/(Y^2 - f_\EC(X),F(X)))\)
%            \;
%            \(Q_1 \gets \pi_p(P)\)
%            ;
%            \(Q_2 \gets \pi_p(Q_1)\)
%            ;
%            \(Q_3 \gets \pi_p(Q_2)\)
%            ;
%            \(Q_4 \gets \pi_p(Q_3)\)
%            ;
%            \(Q_5 \gets [q\bmod\ell]P\)
%            \;
%            \If{\(Q_5 \not= -Q_4\)}{
%                \(Q_6 \gets Q_4 \oplus Q_5\)
%                \;
%                \(t_\ell \gets \textsc{DiscreteLogarithm}(Q_2,Q_6)\)
%                \;
%            }
%            \lElse{
%                \(t_\ell \gets 0\)
%            }
%        }
    }
    \Return{\(\textsc{ChineseRemainderTheorem}(\mathcal{T})\)}
    \;
\end{algorithm}

\section{%%%%%%%%%%%%%%%%%%%%%%%%%%%%%%%%%%%%%%%%%%%%%%%%%%%%%%%%%%%%%%%%%%%%%%%
    Admissible curves
}%%%%%%%%%%%%%%%%%%%%%%%%%%%%%%%%%%%%%%%%%%%%%%%%%%%%%%%%%%%%%%%%%%%%%%%%%%%%%%%
\label{sec:the-endomorphism}

\begin{center}
    \textbf{From now on, \(\bf q = p^2\).}
\end{center}

Recalling Definition~\ref{def:admissible}:
let \(\EC\) be a \(d\)-admissible curve over~\(\FF_{p^2}\),
with separable \(d\)-isogeny
\(\phi: \EC \to \pconj{\EC}\)
(satisfying \(\pconj{\phi} = \epsilon\dualof{\phi}\) with
\(\epsilon=\pm1\)),
and associated endomorphism \(\psi = \pconj{\phi} \circ \pi_p\).
\begin{proposition}
    The associated 
    and Frobenius endomorphisms of \(\EC\)
    are related by
    \begin{equation}
        \label{eq:psi-squared}
        \psi^2 
        % = \pi_p\circ\phi \circ \pi_p\circ\phi
        =
        [\epsilon d] \pi_{p^2}
        \ .
    \end{equation}
    The characteristic polynomial of \(\psi\) is
    \begin{equation}\label{charpoly}
        \chi_{\psi}(T)
        =
        T^2 - rd T + dp
        \ ,
    \end{equation}
    where \(r\) is an integer satisfying
    \begin{equation}
        \label{eqr}
        dr^2
        =
        2p + \epsilon t_\EC
        \ .
    \end{equation}
    In particular,
    \begin{equation}
        \label{eq:psi-relation}
        r \psi = p + \epsilon\pi_{p^2}
        \quad
        \text{in}
        \quad 
        \End(\EC)
        \ .
    \end{equation}
\begin{proof}
    Equation~\eqref{eq:psi-squared}
    holds because
    \(
        \psi^2
        =
        (\pconj{\phi}\pi_p)(\pconj{\phi}\pi_p)
        =
        (\epsilon\dualof{\phi}\phi)(\pconj{\pi_p}\pi_p)
        =
        [\epsilon d]\pi_q
    \).
    The degree of \(\psi\) is \(dp\),
    so \(\psi\) has characteristic polynomial
    \(\chi_{\psi}(T) = T^2 - xT + dp\)
    for some integer \(x\).
    On the other hand, 
    \(\epsilon d\pi_{q}\) has characteristic polynomial
    \(T^2 - \epsilon d t_{\EC} T + d^2p^2\);
    but \(\psi^2 = x\psi - dp\) is a root,
    so
    \(x = rd\)
    where \(r\) satisfies~\eqref{eqr}.
    We then have \(\epsilon dr^2\pi_{q} = (\epsilon\pi_{q} + p)^2\) 
    in \(\ZZ[\pi_{q}]\).
    Comparing with~\eqref{eq:psi-squared},
    we find \(r\psi = \pm(p + \epsilon\pi_{q})\);
    but then \(\chi_{\psi}(\psi) = 0\)
    implies~\eqref{eq:psi-relation}.
\end{proof}
\end{proposition}

Equation~\eqref{eqr} has a number of interesting corollaries.
First, 
\(t_\EC\equiv-\epsilon2p \pmod{d}\),
so we obtain some information on \(t_\EC\) for free.
Second,
\(r\) determines \(t_\EC\),
and hence \(\#\EC(\FF_{p^2})\).
Third,
we have a much smaller bound on \(r\) than on \(t_\EC\):
for \(d\)-admissible curves the Hasse--Weil bound becomes
\begin{equation}
    \label{eqr-bound}
    |r| \le 2\sqrt{p/d} 
    \ .
\end{equation}
This suggests our point-counting strategy,
which is to modify the SEA algorithm to 
compute \(r\) instead of \(t_\EC\),
by considering the action on \(\EC[\ell]\)
of $\psi$ instead of~$\pi_{q}$
and using fewer primes \(\ell\).

We simplify the task by quickly disposing of the supersingular case,
which can be efficiently detected using
Sutherland's algorithm~\cite{Sutherland12},
or slightly faster using a probabilistic algorithm.
\begin{proposition}
    \label{prop:supersingular-r}
    If \(\EC/\FF_{p^2}\) is \(d\)-admissible,
    then it is supersingular
    if and only if \(r = 0\),
    in which case \(t_\EC = -2\epsilon p\)
    and \(\EC(\FF_{p^2}) \cong (\ZZ/(p + \epsilon)\ZZ)^2\).
\begin{proof}
    The curve \(\EC\) is supersingular if and only if \(p\mid t_{\EC}\),
    if and only if \(p\mid r\)
    (by~\eqref{eqr} mod~\(p\) and \(p\nmid d\)),
    if and only if
    \(r = 0\) (by~\eqref{eqr-bound}).
    The group structure follows 
    from~\cite[Th.~1.1]{Zhu00}.
\end{proof}
\end{proposition}

From now on, we will {\bf assume \(\EC\) is ordinary};
so \(\End(\EC)\) is 
an order in the quadratic imaginary field \(\QQ(\pi_q)\), 
and \(\ZZ[\pi_{q}]\) and \(\ZZ[\psi]\) are orders contained in \(\End(\EC)\). 
Looking at~\eqref{charpoly}, 
we see that the discriminants of \(\ZZ[\psi]\) and \(\ZZ[\pi_q]\)
are related by
\[
    \Delta_{\psi} = d(dr^2 - 4p) 
    \qquad 
    \text{and}
    \qquad
    \Delta_{\pi_{q}} = t_\EC^2 - 4p^2 = r^2\Delta_{\psi}
    \ ,
\]
so 
\(|r|\) is the conductor of \(\ZZ[\pi_{q}]\)
in \(\ZZ[\psi]\): that is,
\[
    \ZZ[\pi_{q}]\subset\ZZ[\psi]\subseteq\End(\EC)
    \qquad
    \text{with}
    \qquad
    [\ZZ[\psi]:\ZZ[\pi_{q}]] = |r|
    \ .
\]
Indeed, since \(\EC\) is ordinary,
we have \(r \not= 0\);
so we can rewrite~\eqref{eq:psi-relation} as
\begin{equation}
    \label{eq:psi-in-EndE}
    \psi = \frac{p + \epsilon\pi_{q}}{r}
    \quad 
    \text{in } \End(\EC)
    \ .
\end{equation}
Deuring's theorem on isogeny classes and class groups
(cf.~\cite[\S4]{Schoof87})
can be used to show that the number of
\(\FF_q\)-isomorphism classes of ordinary \(d\)-admissible curves
with a given \(r\) is \(H(\Delta_\psi)\),
where \(H\) is the Kronecker class number.
In particular,
every \(r\) in the interval of~\eqref{eqr-bound} 
occurs for some \(d\)-admissible \(\EC/\FF_{q}\).

In the language of isogeny volcanoes~\cite{FoMo02}: 
if \(\ell\) is a prime dividing \(r\),
then \(\EC\) is somewhere \emph{strictly above} 
the floor of the volcano for \(\ell\);
that is,
all \(\ell\mid r\) are upper-volcanic.

%%%%% S
\section{%%%%%%%%%%%%%%%%%%%%%%%%%%%%%%%%%%%%%%%%%%%%%%%%%%%%%%%%%%%%%%%%%%%%%%%
    Computing the cardinality of admissible curves
}%%%%%%%%%%%%%%%%%%%%%%%%%%%%%%%%%%%%%%%%%%%%%%%%%%%%%%%%%%%%%%%%%%%%%%%%%%%%%%%

Let \(\EC/\FF_q\) be an ordinary \(d\)-admissible curve,
with associated endomorphism \(\psi\);
we want to compute \(\#\EC(\FF_q)\).
Many of the techniques used in the conventional SEA algorithm can be
transposed to working with \(\psi\) instead of \(\pi_q\).
Equations~\eqref{eqr} and~\eqref{eqr-bound}
show that \(t_\EC\) is completely determined by
\(|r|\), which is bounded by \(2\sqrt{p/d}\);
so we can compute $t_\EC$ 
by computing
\[
    r_\ell := r \pmod{\ell}
\]
for \(\ell\) in a collection of small primes \(\mathcal{L}\)
such that
\[
    \prod_{\ell \in \mathcal{L}} \ell > 4\sqrt{p/d}
    \ ,
\]
then recovering \(r\) from the \(r_\ell\) using the CRT.
As a quick comparison, 
using SEA with $\pi_{q}$ to compute \(t_\EC\) directly 
would require
$\prod_{\ell \in \mathcal{L}} \ell > 4 \sqrt{q} = 4 p$.

\begin{proposition}
    \label{prop:AdmissibleTrace-correctness}
    If \(\EC/\FF_{p^2}\)
    is \(d\)-admissible,
    then under the standard heuristic on prime classes,
    Algorithm~\ref{alg:AdmissibleTrace}
    computes \(t_\EC\) in 
    \(\widetilde{O}(\log^3p)\) expected \(\FF_q\)-operations
    (that is \(\widetilde{O}(\log^4p)\) expected bit operations,
    using fast arithmetic).
\begin{proof}
    We compute \(t_\EC\) from \(r\),
    which we recover exactly using the CRT
    from the pairs \((r_\ell,\ell)\) in \(\mathcal{R}\),
    since \(\prod_{(r_\ell,\ell)\in\mathcal{R}}\ell > 4\sqrt{p/d}\).
    Our approach for computing \(r_\ell\) depends
    on which class \(\ell\) falls into;
    we determine the class of \(\ell\)
    in Lines~\ref{line:ATrace:class}, \ref{line:ATrace:is-volcanic},
    and~\ref{line:ATrace:is-Elkies}
    (exactly as in Algorithm~\ref{alg:SEA}).

    If \(\ell\) is volcanic
    (Lines~\ref{line:ATrace:volcanic-start}
    to~\ref{line:ATrace:volcanic-end}),
    then combining \(\ell\mid\Delta_{\pi_q}\) with~\eqref{eqr}
    yields 
    \(r \equiv 0\) or \(\pm2\sqrt{p/d}\pmod{\ell}\);
    in particular, 
    if \(\ell\) is volcanic and \(dp\) is a nonsquare modulo \(\ell\),
    then \(r_\ell = 0\).

    If \(\ell\) is Elkies
    (Lines~\ref{line:ATrace:Elkies-start}
    to~\ref{line:ATrace:Elkies-end}),
    then let 
    \(
        \EC[\ell] 
        = 
        (\ker\isogeny_1)\oplus(\ker\isogeny_2)
    \)
    be the decomposition of the \(\ell\)-torsion
    into eigenspaces for \(\pi_q\).
    Since \(\ell\) is not volcanic 
    we have \(r \not\equiv 0\pmod{\ell}\),
    so~\eqref{eq:psi-in-EndE}
    shows that the \(\ker\isogeny_i\) are also eigenspaces for \(\psi\).
    So let \(\lambda_{\pi}\) and \(\lambda_\psi\) be the eigenvalues of
    \(\pi_q\) and \(\psi\) on \(\ker\isogeny_1\) (say);
    then~\eqref{eq:psi-in-EndE} yields
    \(
        \lambda_\psi 
        \equiv 
        (p + \epsilon\lambda_{\pi})/r
        \pmod{\ell}
    \),
    and then \(\chi_\psi(\lambda_\psi) \equiv 0 \pmod{\ell}\)
    implies
    \(
        r_\ell
        \equiv 
        \frac{\lambda_\psi}{d} + \frac{p}{\lambda_\psi} 
        \pmod{\ell}
    \).
    We can therefore compute \(r_\ell\)
    by computing \(\lambda_\psi\),
    which is the discrete logarithm 
    of \(\psi(P)\) to the base \(P\)
    for a symbolic point \(P\) in \(\ker\isogeny_1\).

    If \(\ell\) is Atkin
    then we skip it completely,
    as in Algorithm~\ref{alg:SEA}
    (but see~\S\ref{sec:complements}).

    In terms of \(\FF_q\)-operations,
    determining the class of \(\ell\) costs
    \(O(\ell^2(\log\ell)^3\log\log\ell + \log q \sfM(\ell))\),
    while computing \(r_\ell\) costs
    \(O(
        (\log p + \log\ell)\sfM(\ell)\log\ell
        + 
        \ell^{(\omega+1)/2}
    )\) 
    if \(\ell\) is volcanic,
    and \(O((\log p + \ell^{1/2})\sfM(\ell) + \ell^{(\omega+1)/2})\)
    if \(\ell\) is Elkies.
    The standard heuristic on prime classes 
    tell us that we will try \(O(\log p)\) primes \(\ell\),
    the largest of which are in \(O(\log p)\);
    so the total complexity is \(\widetilde{O}(\log^3p)\)
    \(\FF_q\)-operations, as claimed.
\end{proof}
\end{proposition}

\begin{algorithm}
    \caption{\textsc{AdmissibleTrace}}
    \label{alg:AdmissibleTrace}
    \SetKwInOut{Input}{Input}\SetKwInOut{Output}{Output}
    \Input{A \(d\)-admissible curve \(\EC/\FF_{p^2}\), where \(p\) is large}
    \Output{The trace of Frobenius of \(\EC\)}
    \(\mathcal{R} \gets \{ \}\)
        \tcp*{\(\mathcal{R}\) will contain the pairs \((r\pmod{\ell},\ell)\)}
    \(M \gets 1\) 
        \tcp*{After each iteration, \(r\) is known modulo \(M\)}
    \(\ell \gets 1\)
    \;
    \While{\(M \le 4\sqrt{p/d}\)}{
        \lRepeat{
            \(\ell\nmid d\)
        }{
            \(\ell \gets \textsc{NextPrime}(\ell)\) 
        }
        \(J \gets \textsc{Gcd}(x^{p^2}-x,\textsc{EvaluatedModularPolynomial}(\ell,\EC))\)
        \; \label{line:ATrace:class}
        \If(\tcp*[f]{\(\ell\) is volcanic}){\(\deg J = 1\) or \(\ell + 1\)}{
            \label{line:ATrace:is-volcanic}
            \If{\(dp\) has a square root \(s\) modulo \(\ell\)}{
 	         \(F \gets
            	 \textsc{KernelPolynomial}(\ell,\EC,\textsc{OneRoot}(J))\)
            	 \; \label{line:ATrace:volcanic-start}
            	 \(P \gets (X,Y)\in \EC(\FF_q[X,Y]/(Y^2 - f_\EC(X),F(X)))\)
            	 \;
                \(Q_1 \gets \pi_p(P)\) 
                ;
                \(Q_2 \gets \pconj{\phi}(Q_1)\)
                ;
                \(Q_3 \gets [s]P\)
                \;
                \(
                    r_\ell 
                    \gets 
                    \begin{cases} 
                        2s/d \pmod{\ell}  & \text{if } Q_3 = Q_2 \ ; \\
                        -2s/d \pmod{\ell} & \text{if } Q_3 = -Q_2 \\
                        0                 & \text{otherwise} \ ;
                    \end{cases}
                \)
                \;
            }
            \lElse{
                \(r_\ell \gets 0\)
            }
            \(\mathcal{R}\gets\mathcal{R}\cup\{(r_\ell,\ell)\}\)
            ;
            \(M \gets \ell M\)
            \; \label{line:ATrace:volcanic-end}
        }
        \ElseIf(\tcp*[f]{\(\ell\) is Elkies}){\(\deg J = 2\)}{
            \label{line:ATrace:is-Elkies}
            \(F \gets
            \textsc{KernelPolynomial}(\ell,\EC,\textsc{OneRoot}(J))\)
            \; \label{line:ATrace:Elkies-start}
            \(P \gets (X,Y)\in \EC(\FF_q[X,Y]/(Y^2 - f_\EC(X),F(X)))\)
            \;
            \(Q_1 \gets \pi_p(P)\)
            ;
            \(Q_2 \gets \pconj{\phi}(Q_1)\)
            \;
            \(r_\ell \gets \lambda/d+p/\lambda \pmod{\ell}\)
            \textbf{where}
            \(\lambda = \textsc{DiscreteLogarithm}(P,Q_2)\)
            \;
            \(\mathcal{R}\gets\mathcal{R}\cup\{(r_\ell,\ell)\}\)
            ;
            \(M \gets \ell M\)
            \; \label{line:ATrace:Elkies-end}
        }
%        \Else(\tcp*[h]{\(\ell\) is Atkin}){
%            \textbf{pass}
%            \;
%            \(F \gets \textsc{DivisionPolynomial}(\ell,\EC)\)
%            \;
%            \(P \gets (X,Y)\) 
%            in \(\EC(\FF_q[X,Y]/(Y^2 - f_\EC(X),F(X)))\)
%            \;
%            \(Q_1 \gets \pi_p(P)\)
%            ;
%            \(Q_2 \gets \pi_p(Q_1)\)
%            ;
%            \(Q_3 \gets [p\bmod\ell]P\)
%            \;
%            \If{\(Q_3 \not= -\epsilon Q_2\)}{
%                \(Q_4 \gets \pconj{\phi}(Q_1)\)
%                ;
%                \(Q_5 \gets \epsilon Q_2 \oplus Q_3\)
%                \;
%                \(r_\ell \gets \textsc{DiscreteLogarithm}(Q_4,Q_5)\)
%                \;
%            }
%            \lElse{
%                \(r_\ell \gets 0\)
%            }
%        }
    }
    \Return{ \(\epsilon(dr^2 - 2p)\) 
        where \(r = \textsc{ChineseRemainderTheorem}(\mathcal{R})\)
    }
    \;
\end{algorithm}

\begin{remark}
    Suppose $\ell \mid d$ and \(\ell \not=2\). 
    Equation~\eqref{eqr} tells us that 
    $t_\EC \equiv 2 \epsilon p \pmod{\ell}$;
    so $\ell \mid \Delta_{\pi_q}$,
    and $\ell$ is volcanic. 
    Moreover, since $\Delta_\psi = d (d r^2 - 4p)$, 
    we can deduce that $\ell \mid\mid\Delta_\psi$. % if \(\ell \not= 2\).
    Note also that
    \(\End(\EC)\cong\End(\pconj{\EC})\),
    so the \(\ell\)-isogeny factoring \(\phi\) is horizontal;
    this implies that \(\End(\EC)\) is \(\ell\)-maximal.
    Combined with the above, we see that \(\ZZ[\psi]\) is
    \(\ell\)-maximal in \(\QQ(\pi_q)\).
    In particular, if \(\ell\) is upper-volcanic then \(\ell\mid r\)
    (and \((0,\ell)\) can be added to \(\mathcal{R}\) in
    Algorithm~\ref{alg:AdmissibleTrace}).
\end{remark}

\section{%%%%%%%%%%%%%%%%%%%%%%%%%%%%%%%%%%%%%%%%%%%%%%%%%%%%%%%%%%%%%%%%%%%%%%%
    Complements
}%%%%%%%%%%%%%%%%%%%%%%%%%%%%%%%%%%%%%%%%%%%%%%%%%%%%%%%%%%%%%%%%%%%%%%%%%%%%%%%
\label{sec:complements}

Schoof's original algorithm may be generalized from prime \(\ell\) 
to small prime powers in a very simple way.
Going further, we may use isogeny cycles to
compute eigenspaces of \(\pi_{q}\) and \(\psi\) on $\EC[\ell^n]$
for Elkies \(\ell\):
the methods developed 
for~\(\pi_{q}\) 
in~\cite{CoMo94} 
and~\cite{GaMo06}
generalize to \(\psi\) without any difficulty. 
Once we have recovered 
\[
    \psi(P)
    = 
    [k_n]P
    \quad
    \text{for}
    \quad
    P = (X,Y)
    \in
    \EC\left(\FF_q[X,Y]/(Y^2-f_\EC(X),F_{\ell^n}(X))\right)
    ,
\]
we have $k_{n+1} = k_n + \tau \ell^n$ for $0\leq \tau < \ell$, 
and we need to test
\[
    \psi(P) - [k_n]P
    = 
    [\tau]([\ell^n]P)
    \quad \text{in} \quad
    \EC\left(\FF_q[X,Y]/(Y^2-f_\EC(X),F_{\ell^{n+1}}(X)\right)
\]
(here \(F_{\ell^n}\) and \(F_{\ell^{n+1}}\) 
are factors of \(\Psi_{\ell^n}\) and \(\Psi_{\ell^{n+1}}\)
that are minimal polynomials for \(\ell^n\) and \(\ell^{n-1}\)-torsion
points).

We may extend Algorithms~\ref{alg:SEA}
and~\ref{alg:AdmissibleTrace}
to use Atkin primes.
If \(\ell\) is Atkin,
then $\pi_q$ and $\psi$ have no rational eigenspaces in \(\EC[\ell]\);
but we may still
compute \(t_\ell\) and \(r_\ell\)
by working on the full \(\ell\)-torsion,
as in Schoof's original algorithm.
If \(P\) is a symbolic point of \(\EC[\ell]\)
then 
\((\chi_{\pi_q} \bmod{\ell})(P) = 0\),
so in Algorithm~\ref{alg:SEA},
\(t_\ell\) is the discrete logarithm
of \(\pi_q(\pi_q(P))+[q\bmod\ell]P\)
to the base \(\pi_q(P)\);
similarly,
in Algorithm~\ref{alg:AdmissibleTrace},
\(r_\ell\) is the discrete logarithm 
of \(\epsilon\pi_{q}(P) + [p \bmod\ell](P)\)
to the base \(\psi(P)\)
(here we use \(\psi^2 - dr\psi + [dp] = d(\epsilon\pi_q - r\psi +
[p]) = 0\) and \(\ell \nmid d\)).
The kernel polynomial defining \(\EC[\ell]\) is \(\Psi_\ell\),
which we can compute using standard recurrences
involving the coefficients of \(f_\EC\)
(using the method of~\cite{Cheng03b}, for example)
in \(O(\sfM(\ell^2)\log \ell)\)
\(\FF_q\)-operations.
But \(\Psi_\ell\) has degree \((\ell^2-1)/2\),
so
computing \(t_\ell\) resp. \(r_\ell\)
costs \(O((\log q)\sfM(\ell^2))\) resp. \(O((\log p)\sfM(\ell^2))\) 
\(\FF_q\)-operations;
for that cost, we would gain much more information 
by using a larger Elkies prime instead. Alternatively, we can use
Atkin's initial ideas using the splitting degree of $\Phi_\ell(X, j(\EC))$
to determine a list of potential $t_\ell$ to be used in a
tricky match and sort algorithm, or the more advanced algorithm
of~\cite{JoLe01}.
In our setting, we could use \eqref{eqr} to transform
the list of $t_\ell$'s to build a list of $r_\ell$'s (on average, this
does not increase the size of the lists too much).

Finally, we mention the use of the baby-step
giant-step approach to speed up the final computations.
If $P \in \EC(\FF_q)$, 
then $\chi_\psi(P) = 0$ becomes
\(
    [\epsilon d + d p] P = [r d] \psi(P)
\),
so 
\([p+\epsilon] (Q) = [r] \psi(Q)\)
with $Q = [d] P$ (if $Q = O_\EC$, then another $P$ should be used). 
Suppose we stop the loop of Algorithm~\ref{alg:AdmissibleTrace} early;
then $r$ is known modulo $M$. 
Writing $r = r_0 + s M$ with $|s| \leq 2 \sqrt{p/d}/M$, 
we can find $s$ 
by solving $[p+\epsilon-r_0] Q = [s] ([M] \psi(Q))$
for a sufficiently general choice of \(Q\) in \(\EC(\FF_q)\);
this is a classical discrete logarithm problem with \(\FF_q\)-points, 
but in a smaller search space than the whole of \(\EC(\FF_q)\). 
The optimal threshold for \(M\) is best determined through experiments.

%%%%% S
\section{%%%%%%%%%%%%%%%%%%%%%%%%%%%%%%%%%%%%%%%%%%%%%%%%%%%%%%%%%%%%%%%%%%%%%%%
    Comparison of Algorithms~\ref{alg:SEA} and~\ref{alg:AdmissibleTrace}
}%%%%%%%%%%%%%%%%%%%%%%%%%%%%%%%%%%%%%%%%%%%%%%%%%%%%%%%%%%%%%%%%%%%%%%%%%%%%%%%
\label{sec:complexity}

Let us compare the cost of computing \(t_\EC\)
with Algorithms~\ref{alg:SEA} and~\ref{alg:AdmissibleTrace}
when \(\EC\) is \(d\)-admissible.
For simplicity,
we will suppose that Algorithm~\ref{alg:SEA} also avoids the primes
dividing \(d\)
(these are very few and very small,
so they do not contribute asymptotically
or practically to the comparison).

The first clear difference between 
the algorithms
is the number and size of primes \(\ell\) used:
Algorithm~\ref{alg:AdmissibleTrace} essentially uses the smaller
half of the set of primes used by Algorithm~\ref{alg:SEA}.
The largest primes in each set still have roughly
the same size, \(O(\log p)\),
so asymptotically this makes no difference---but using half the number of
primes, and the smaller half at that, represents an important
improvement in practice.

Now consider the cost of computing \(t_\ell\)
(as in Algorithm~\ref{alg:SEA})
or \(r_\ell\)
(as in Algorithm~\ref{alg:AdmissibleTrace})
for the same \(\ell\).
The costs of determining the class of \(\ell\)
and the calls to
\(\textsc{KernelPolynomial}\)
are identical,
and the calls to
\(\textsc{DiscreteLogarithm}\) are equivalent.
The only real difference is 
in how each algorithm computes the relations 
used to
determine \(t_\ell\) and \(r_\ell\).
\begin{itemize}
    \item \underline{if \(\ell\) is Elkies},
        then
        Algorithm~\ref{alg:SEA} uses 
        \(2\times \pi_p\)
        while
        Algorithm~\ref{alg:AdmissibleTrace} uses
        \( 1\times \pi_p + 1\times \pconj{\phi} \).
%    \item \underline{if \(\ell\) is Atkin},
%        then
%        Algorithm~\ref{alg:SEA} 
%        uses \(4\times \pi_p + 1\times[q\bmod\ell] + 1\times\oplus\),
%        while
%        Algorithm~\ref{alg:AdmissibleTrace}
%        uses \(2\times\pi_p + 1\times\pconj{\phi} +
%        1\times[p\bmod\ell] + 1\times\oplus\).
    \item \underline{if \(\ell\) is volcanic},
        then (in the worst cases)
        Algorithm~\ref{alg:SEA} uses \(2\times\pi_p + 1\times[s\bmod\ell]\),
        while Algorithm~\ref{alg:AdmissibleTrace}
        uses \(1\times \pi_p + 1\times \pconj{\phi} + 1\times[s\bmod\ell]\).
\end{itemize}
In each case, 
the asymptotic costs are the same;
but if \(d \ll \log p\),
then the costs are dominated by computations of \(\pi_p\)
on \(\subgrp{P}\) (for the same \(P\)).
The crucial practical difference 
is that for each class of prime,
Algorithm~\ref{alg:AdmissibleTrace}
exchanges half of the computations of \(\pi_p\) required by
Algorithm~\ref{alg:SEA}
for one computation of \(\pconj{\phi}\),
which has a very small cost when \(d \ll \log p\).
Hence, for any given prime~\(\ell\), 
Algorithm~\ref{alg:AdmissibleTrace}
should compute \(r_\ell\)
twice as quickly as Algorithm~\ref{alg:SEA}
computes \(t_\ell\).

By our complexity analysis, we see that the largest $\ell$ is $O(\log
p)$ instead of $O(\log q)$, and we use the the smaller half of them,
we expect a real speedup of a factor of four. This is confirmed by our
experimental results in \S\ref{sec:implementation} below.

%%%%% S
\section{%%%%%%%%%%%%%%%%%%%%%%%%%%%%%%%%%%%%%%%%%%%%%%%%%%%%%%%%%%%%%%%%%%%%%%%
    $\Q$-curves and other sources of admissible curves
}%%%%%%%%%%%%%%%%%%%%%%%%%%%%%%%%%%%%%%%%%%%%%%%%%%%%%%%%%%%%%%%%%%%%%%%%%%%%%%%
\label{sec:QQ-curves}

Admissible curves appear naturally as reductions of quadratic
\(\QQ\)-curves modulo inert primes (cf.~\cite[\S3]{Smith15}).
As such, we can construct parametrized families of admissible curves
over any~\(\FF_{p^2}\).

\begin{definition}
    \label{def:QQ-curve}
    A \emph{quadratic \(\QQ\)-curve of degree \(d\)}
    is an elliptic curve \(\ECK\) 
    without complex multiplication,
    defined over a quadratic field \(\QQ(\sqrt{\Delta})\),
    such that there exists an isogeny of degree \(d\)
    from \(\ECK\) to its Galois conjugate \(\conj[\tau]{\ECK}\),
    where \(\tau\) is the conjugation of \(\QQ(\sqrt{\Delta})\) over \(\QQ\).
\end{definition}

\begin{proposition}
    \label{prop:QQ-curve-reductions-are-admissible}
    Let \(\ECK/\QQ(\sqrt{\Delta})\) be a quadratic \(\QQ\)-curve of degree
    \(d\).
    If \(p\nmid d\) is a prime of good reduction for \(\ECK\)
    that is inert in \(\QQ(\sqrt{\Delta})\), 
    then
    the reduction of \(\ECK\) modulo \(p\) is \(d\)-admissible.
\begin{proof}
    Gonz\'alez shows that
    a \(d\)-isogeny \(\phiK:\ECK\to\conj[\tau]{\ECK}\) 
    must be defined over \(\QQ(\sqrt{\Delta},\sqrt{\pm d})\)
    (see~\cite[\S3]{Gonzalez01});
    so if we extend \(\tau\) 
    to the involution of \(\QQ(\sqrt{\Delta},\sqrt{\pm d})\)
    that acts trivially on \(\QQ(\sqrt{\pm d})\)
    if and only if \(\sqrt{\pm d}\) is in \(\FF_p\),
    then \(\phiK\)
    reduces modulo \(p\) to a \(d\)-isogeny \(\phi: \EC \to \pconj{\EC}\)
    over \(\FF_{p^2}\),
    and \(\conj[\tau]{\phiK}\) reduces to \(\pconj{\phi}\).
    Observe that \(\conj[\tau]{\phiK}\phiK\) is an endomorphism of
    \(\ECK\) of degree \(d^2\).
    Since \(\ECK\) does not have complex multiplication,
    its only endomorphisms of degree \(d^2\)
    are \([\pm d]\);
    hence 
    \(\conj[\tau]{\phiK} = \epsilon\dualof{\phiK}\)
    with \(\epsilon = \pm1\).
    Reducing modulo \(p\) we have \(\pconj{\phi} = \epsilon\dualof{\phi}\),
    so \(\EC\) is \(d\)-admissible.
\end{proof}
\end{proposition}

We emphasize that if a \(d\)-admissible curve \(\EC\) is the reduction
of a quadratic \(\QQ\)-curve \(\ECK\),
then the associated endomorphism on \(\EC\) is \emph{not} the reduction
of any endomorphism on \(\ECK\).
Indeed, \(\ECK\) has no non-integer endomorphisms by definition.

\begin{example}
    \label{ex:Hasegawa-2}
    Fix any prime \(p > 3\);
    the following construction
    (carried much further in~\cite{Smith13}
    and~\cite{Smith15})
    yields a 1-parameter family of
    \(2\)-admissible curves over \(\FF_{p^2}\).
    Let \(\Delta\) be a squarefree integer that is not a square modulo \(p\)
    (so \(p\) is inert in \(\QQ(\sqrt{\Delta})\)),
    let \(\tau\) be the involution of \(\QQ(\sqrt{\Delta},\sqrt{-2})\)
    that restricts to \(\sigma\) modulo \(p\),
    and let \(s\) be a free parameter taking values in \(\QQ\).
    The family of curves over \(\QQ(\sqrt{\Delta})\)
    defined by
    \(
        \ECK : 
        y^2 
        = 
        x^3 - 6(5-3s\sqrt{\Delta})x + 8(7 - 9s\sqrt{\Delta}) 
    \)
    is equipped with a 2-isogeny \(\phiK: \ECK \to \conj[\tau]{\ECK}\)
    over \(\QQ(\sqrt{\Delta},\sqrt{-2})\)
    with kernel polynomial \(D(x) = x - 4\)
    (see~\cite[Prop 3.3]{Hasegawa97}).
    Computing \(\dualof{\phiK}\) and \(\conj[\tau]{\phiK}\),
    we find that
    \(\conj[\tau]{\phiK} = \epsilon\dualof{\phiK}\),
    where \(\epsilon = 1\) if \(p \equiv 5, 7 \pmod{8}\)
    and \(\epsilon = -1\) if \(p \equiv 1, 3 \pmod{8}\).
    Reducing everything modulo \(p\), 
    as in the proof of Prop.~\ref{prop:QQ-curve-reductions-are-admissible},
    we obtain a family of curves
    \[
        \EC : 
        y^2 
        = 
        x^3 - 6(5-3s\sqrt{\Delta})x + 8(7 - 9s\sqrt{\Delta}) 
        \quad
        \text{over}
        \quad
        \FF_{p^2} = \FF_p(\sqrt{\Delta})
    \]
    with the parameter \(s\) taking values in \(\FF_p\),
    equipped with a 2-isogeny \(\phi: \EC \to \pconj{\EC}\)
    over \(\FF_{p^2}\).
    Composing \(\pi_p\) with \(\pconj{\phi}\)
    yields the associated endomorphism \(\psi\) of \(\EC\),
    defined by
    \[
        \psi:
        (x,y) 
        \longmapsto 
        \left(
            \frac{x^p(x^p - 4) + 18(1 - s\sqrt{\Delta})}{-2(x^p-4)}
            \ ,
            \frac{y^p}{\sqrt{-2}^p}\left(
                \frac{(x^p-4)^2 - 18(1 - s\sqrt{\Delta})}{-2(x^p-4)^2}
            \right)
        \right)
        \ .
    \]
\end{example}

Since the definition of admissible curves involves only isogenies over
\(\FF_{p^2}\), we would expect a characterization of admissible curves
over a given \(\FF_{p^2}\) in terms of modular polynomials.

\begin{proposition}
    If \(\EC\) is an ordinary elliptic curve
    over \(\FF_q = \FF_{p^2}\)
    such that \(j(\EC)\) is 
    a simple root of \(\Phi_d(x,x^p)\) in
    \(\FF_{q}\setminus\{0,1728\}\)
    (so in particular,
    \(\Aut_{\FFbar_{q}}(\EC) = \{[\pm1]\}\)),
    then \(\EC\) is \(d\)-admissible.
\begin{proof}
    If \(j(\EC)\) is a simple root of \(\Phi_d(x,x^p)\) in \(\FF_{q}\),
    then up to \(\FFbar_{q}\)-isomorphism there is a unique
    \(d\)-isogeny
    \(\phi: \EC \to \pconj{\EC}\).
    If \(\phi\) were not defined over \(\FF_q\),
    then the endomorphism \(\pconj{\pi_p}\phi\)
    would not be defined over \(\FF_{q}\),
    hence not commute with \(\pi_{q}\),
    contradicting non-supersingularity.
    For \(d\)-admissibility,
    it remains to show that \(\pconj{\phi} = \epsilon\dualof{\phi}\)
    with \(\epsilon = \pm1\). % FIXME: can we say \epsilon right now?
    But if this were not the case,
    then \(\dualof{(\pconj{\phi})}\) would be a second \(d\)-isogeny
    \(\EC\to\pconj{\EC}\),
    %defined over \(\FF_{q}\),
    not isomorphic to \(\phi\)
    (since \(\Aut_{\FFbar_{q}}(\EC) = \{[\pm1]\}\));
    that is, \(j(\EC)\) would be (at least) a double root of
    \(\Phi_d(x,x^p)\).
\end{proof}
\end{proposition}

\begin{example}
    Multiple roots of \(\Phi_d(x,x^p)\)
    may not yield \(d\)-admissible curves.
    Consider
    the ordinary curve
    \(\EC: y^2 = x^3 + (38 + 53i)x + 27 - 3i\)
    over \(\FF_{q} = \FF_{103}(i)\)
    where \(i^2 = -1\):
    then \(j(\EC) = 35 + 5i\)
    is a double root of \(\Phi_3(x,x^{103})\).
    Indeed, we have a pair of non-isomorphic \(3\)-isogenies
    \(\phi_1 : \EC \to \pconj{\EC}\)
    and \(\phi_2 : \EC \to \pconj{\EC}\),
    with kernel polynomials \(x + 1 + 39i\)
    and \(x - 4 + 32i\), respectively;
    but
    \(\pconj{\phi_1} = \pm\dualof{\phi_2}\)
    and
    \(\pconj{\phi_2} = \pm\dualof{\phi_1}\),
    so \(\EC\) is not \(3\)-admissible.
\end{example}

%%%%% S
\section{%%%%%%%%%%%%%%%%%%%%%%%%%%%%%%%%%%%%%%%%%%%%%%%%%%%%%%%%%%%%%%%%%%%%%%%
    Generating cryptographically strong curves
}%%%%%%%%%%%%%%%%%%%%%%%%%%%%%%%%%%%%%%%%%%%%%%%%%%%%%%%%%%%%%%%%%%%%%%%%%%%%%%%

One of the important motivations for developing our algorithm
was the generation of cryptographically strong curves.
Indeed, 
the curves proposed for cryptographic applications 
in~\cite{Smith13} and~\cite{Smith15}, 
and which were subsequently used in fast, compact Diffie--Hellman key
exchange software~\cite{CoHiSm14},
are admissible.
These curves were designed to offer accelerated scalar multiplication
(using the associated endomorphism) over fast finite fields,
without obstructing twist-security;
but when generating twist-secure curves
at and above the 128-bit security level, 
we can expect to try hundreds of thousands of curves 
before finding a suitable one.
In this context of counting many curves, 
practical speedups become very important.

For cryptographic applications based on the hardness of the discrete
logarithm problem,
the minimum requirement for a ``secure'' curve \(\EC/\FF_{p^2}\)
is that \(\#\EC(\FF_{p^2}) = c\cdot n\),
where \(n\) is prime and \(c\) is tiny
(traditionally, we want \(c = 1\); more modern software using
Montgomery and Edwards models requires \(c = 2\) or \(4\)).
For some applications we further require ``twist-security'':
that is, the quadratic twist \(\EC'\) should 
satisfy \(\#\EC'(\FF_{p^2}) = c'\cdot n'\),
where \(n'\) is prime and \(c'\) is tiny.

To find a secure or twist-secure curve over \(\FF_{p^2}\)
we typically fix a prime \(p\) 
of bitlength around the required security parameter,
then test a series of curves over \(\FF_{p^2}\),
computing their orders
until we find a curve with the right structure.
Equation~\eqref{eqr} implies
\[
    \#\EC(\FF_{p^2})
    =
    (p + \epsilon)^2 - \epsilon dr^2
    \quad
    \text{and}
    \quad
    \#\EC'(\FF_{p^2})
    =
    (p-\epsilon)^2 + \epsilon dr^2
    \ .
\]
This places some immediate constraints on the combinations of \(d\),
\(p\), and \(\epsilon\) that can yield suitable curves.
For example,
\(\#\EC(\FF_{p^2}) \equiv (p + \epsilon)^2 \pmod{d}\),
so \(d\mid\#\EC(\FF_{p^2})\)
if and only if 
\(p \equiv -\epsilon \pmod{d}\);
such \(p\) should be avoided unless we can accept \(d\mid c\).
Similarly, if twist-security prohibits \(d\mid c'\)
then we should must avoid \(p \equiv \epsilon\pmod{d}\).
Clearly if \(\EC\) is \(2\)-admissible, then it must have a rational point of
order 2, so we cannot do better than having \(c = c' = 2\).
Similarly, \(3\)-admissible curves must have either \(3\mid c\) or \(3\mid c'\).

Extensive computations done for $d = 2$ and $3$ over a range of primes
revealed densities of twist-secure \(d\)-admissible curves (modulo the
constraints above)
similar to the densities of twist-secure general elliptic curves over the same
fields.

With Algorithm~\ref{alg:SEA},
we can speed up the search for secure curves
by checking whether 
\(t_\ell \equiv p^2 + 1 \pmod{\ell}\)
for each \(\ell\);
if so, 
then \(\ell\mid\#\EC(\FF_{p^2})\),
so we can abort the computation 
and move on to the next candidate curve~\cite{Lercier97}.
Similarly, if \(t_\ell \equiv -(p^2 + 1)\pmod{\ell}\)
then \(\ell\mid\#\EC'(\FF_{p^2})\).

With Algorithm~\ref{alg:AdmissibleTrace},
if \(\ell\) divides \(\#\EC(\FF_{p^2})\)
then \((p + \epsilon)^2 \equiv \epsilon d r^2 \pmod{\ell}\),
so
\(\ell\) cannot divide \(\#\EC(\FF_{p^2})\)
unless \(\epsilon d\) is a square mod \(\ell\);
and if \(\epsilon d\) is a square mod \(\ell\),
then we should abort 
if \(r_\ell \equiv \pm (p + \epsilon)/\sqrt{\epsilon d} \pmod{\ell}\).
In fact,
if \(r_\ell \equiv 0\) and \(p + \epsilon \equiv 0 \pmod{\ell}\),
then the nondegeneracy of the \(\ell\)-Weil pairing
%(which maps into the order-\(\ell\) subgroup of \(\FF_{p^2}^\times\)
%in this case)
implies that \(\EC[\ell](\FF_{p^2}) \cong (\ZZ/\ell\ZZ)^2\).
Replacing \(\epsilon\) with \(-\epsilon\)
yields analogous results for the twist \(\EC'\).

We note also that there may be an advantage in generating curves using
the parameter $r$ and not $t_\EC$. We could
force some value of $\ell$ to divide $r$ by rejecting curves $\EC$ for
which $\Phi_\ell(X, j(\EC))$ does not have $1$ or $\ell+1$ roots. This
has no impact on $t_\EC$, and we already know $r\pmod\ell$. We just need
to hope that such curves are as secure as general \(d\)-admissible curves.

%%%%% S
\section{Implementation and experiments}
\label{sec:implementation}

We implemented the new algorithm on top of our implementation of SEA,
realized in C++ using NTL 9.6.4 (with \texttt{gcc} 4.9.2).
The timings below (in seconds)
are for an Intel Xeon platform (E5520 CPU at 2.27GHz).
We define two primes (of 128 and 255 bits), 
derived from the decimal expansion of \(\pi\):
\begin{align*}
    p_{128} & :=\ \scriptstyle 314159265358979323846264338327950288459\ ,
    \\
    p_{255} & :=\ \scriptstyle
    31415926535897932384626433832795028841971693993751058209749445923078164062963
    \ .
\end{align*}

First, we compare the
straightforward computation of $X^q \bmod \Phi_\ell$ to a modular
composition over \(\FF_{p^2}\) with \(p = p_{128}\) and \(p_{255}\),
for two choices of $\ell$:
$$\begin{array}{|c||c|c||c|}\hline
\multicolumn{4}{|c|}{p_{128}} \\ \hline
\ell & X^p \bmod \Phi_\ell& X^p \circ X^p & X^q \\ \hline
101 & 0.23 & 0.04 & 0.47 \\
173 & 0.43 & 0.11 & 0.88 \\
\hline
\end{array}\quad
\begin{array}{|c||c|c||c|}\hline
\multicolumn{4}{|c|}{p_{255}} \\ \hline
\ell & X^p \bmod \Phi_\ell& X^p \circ X^p & X^q \\ \hline
101 & 0.69 & 0.07 & 1.40 \\
173 & 1.38 & 0.18 & 2.80 \\
%% 211 & 1.50 & 0.22 & 3.02 \\ %% not really useful
\hline
\end{array}$$
Then we ran our program 
on curves from the %Hasegawa
family of Example \ref{ex:Hasegawa-2}, 
for each $1\leq s\leq 100$.
This gave the following average values:
$$\begin{array}{|c||c|c|}\hline
\multicolumn{3}{|c|}{p_{128}} \\ \hline
& \text{Algorithm~\ref{alg:SEA}} &
\text{Algorithm~\ref{alg:AdmissibleTrace}} \\ \hline
\text{max. }\ell & 164 & 62 \\
X^q \text{ time} & 9.11 & 2.62 \\
\text{Total time} & 20.11 & 4.1 \\
\hline
\end{array}\quad
\begin{array}{|c||c|c|}\hline
\multicolumn{3}{|c|}{p_{255}} \\ \hline
& \text{Algorithm~\ref{alg:SEA}} &
\text{Algorithm~\ref{alg:AdmissibleTrace}} \\ \hline
\text{max. } \ell & 352 & 160.76 \\
X^q \text{ time} & 89.73 & 22.55 \\
\text{Total time} & 171.95 & 39.16 \\
\hline
\end{array}$$
Finally, we searched for twist-secure curves with small values of
the parameter $s$. For instance, with $p = p_{128}$ and $s=113$, 
we get a curve of cardinality $2 p'$, whose twist has cardinality $6 p''$; 
with $p = p_{255}$, taking $s=269$ yields a pair of curves 
each with cardinality two times a prime.

\medskip
\noindent
{\bf Acknowledgments.} We thank A.~Sutherland for pointing out an error in
the complexity analysis of the SEA algorithm.

\appendix
\section{%%%%%%%%%%%%%%%%%%%%%%%%%%%%%%%%%%%%%%%%%%%%%%%%%%%%%%%%%%%%%%%%%%%%%%%
    Detailed complexity of basic computations
}%%%%%%%%%%%%%%%%%%%%%%%%%%%%%%%%%%%%%%%%%%%%%%%%%%%%%%%%%%%%%%%%%%%%%%%%%%%%%%%
\label{sec:technical-lemmas}

Let $F(X)$ be a degree $e$ polynomial with coefficients in
$\FF_q[X]$. 
We define \(G\) and \(H\) to be the polynomials of degree \(< e\) 
such that
\(H \equiv X^p \pmod{F}\)
and
\begin{equation}\label{eqY}
    Y^p \equiv Y G(X) 
    \text{ with } 
    G(X) \equiv f_\EC^{(p-1)/2}(X) \bmod F(X)
    \ .
\end{equation}

\subsection{Computing $X^q \bmod F$}

The first step in factoring $F$ is to compute $X^q \bmod
F$. When $q = p^n$ for some prime $p$, 
we may start by computing \(H\) and then proceed with modular
composition.

If $R(X) = \sum_{i=0}^{e-1} r_i X^i$ with $r_i \in \FF_q$, 
then $\pconj{R}(X) = \sum_{i=0}^{e-1} r_i^p X^i$ 
satisfies $R^p \bmod F = \pconj{R} \circ X^p \bmod F$.
We assume that the cost of computing all the $r_i^p$ is negligible
(as it is with a suitable choice of basis for \(\FF_q/\FF_p\):
if $\FF_{p^2} = \FF_p(\sqrt{\Delta})$, then
$(a+b\sqrt{\Delta})^p = a-b\sqrt{\Delta}$ for all \(a\) and \(b\) in
\(\FF_p\)).
For our purposes, the computation of $X^{p^2}$ 
computes \(H(X)\) and \(X^{p^2} = \pconj{H}\circ H \bmod F\),
which costs $O((\log p) \sfM(e) + \mathcal{C}(e))$ instead of
$O((\log q) \sfM(e))$, which is larger provided that $2 e \leq (\log
p)^2$. 
When $q = p^n$ with $n > 2$, similar savings can be obtained.

\subsection{Proof of Lemma~\ref{lemma:costs}}

Let \(F = D\), or any factor of \(D\) 
(as in the extensions of the algorithm mentioned
in~\S\ref{sec:complements}).

For (i), 
the obvious way is to compute $H$, then $G$,
in $O((\log p){\sfM}(e))$
\(\FF_q\)-operations.
Alternatively, we can adapt the methods of~\cite{GaMo06}: 
first compute \(G\) 
in $O((\log p) {\sfM}(e))$ operations. 
Consider the polynomial $P(W) = W^3 + A^p W +
B^p - (X^3+A X+B) G(X)^2$. 
Then $X^p \bmod F$ is a root of 
both $\pconj{F}(W)$ and $P(W)$
in $\FF_q[X]/(F(X))$,
so $W - H(X) \mid g = \gcd(P(W), \pconj{F}(W))$. 
Very
generally, $g = W - H(X)$. The main cost is that of
reducing $\pconj{F}(W)$ modulo $P(W)$, which is $O(e \sfM(e))$. 
This can be reduced to $\mathcal{C}_3(e)$
or even $O((\log \ell) \sfM(e))$
if $F$ divides $\Psi_\ell$.

For (ii):
we can compute 
\( 
    \pi_p(Q) 
    = 
    (Q_x^p, Y^p Q_y^p)
    = 
    (\pconj{Q_x}\circ H \bmod F, Y G (\pconj{Q}_y\circ H) \bmod F)
\)
in $\mathcal{C}_2(e)$ \(\FF_q\)-operations.
This also applies for computing 
$\pi_q(P) = (X^{p^2}, Y^{p^2}) = \pi_p (H, Y G)$.

For (iv):
suppose $\phi = (N/D, M/D^2)$ with
\(\deg N = \deg M = 2\) and \(\deg D = 1\). 
We compute
\( N \circ Q_x \bmod F \),
\( M \circ Q_x \bmod F \),
and
\( D \circ Q_x \bmod F \)
followed by some multiplications, keeping numerators and denominators.
We only need a few modular multiplications, for a cost of $O({\sfM}(e))$.

For (v), 
we have $\phi = (N/D^2, M/D^3)$ with
$\mathrm{deg}(N) = \mathrm{deg}(M) = d$, and $\mathrm{deg}(D) = (d-1)/2$.
First we reduce \(N\), \(M\), and \(D\) modulo \(F\) (if necessary),
at a cost of \(O(M(d))\).  
We then compute
\( N \circ Q_x \bmod F \),
\( M \circ Q_x \bmod F \),
and
\( D \circ Q_x \bmod F \)
followed by some multiplications, keeping numerators and denominators.  
The dominating cost is bounded by $O({\sfM}(d) + \mathcal{C}_3(e))$.

\bibliographystyle{plain}
\bibliography{qcsea}

\affiliationone{
F. Morain\\
\'Ecole Polytechnique/LIX\\
and Centre national de la recherche scientifique (CNRS)\\
and Institut national de recherche en informatique et en
automatique (INRIA)\\
France
   \email{morain@lix.polytechnique.fr}}
\affiliationtwo{
C. Scribot\\
Minist\`ere de l'\'Education Nationale \\
France}
\affiliationthree{
B.~Smith\\
Institut national de recherche en informatique et en
automatique (INRIA)\\
and \'Ecole Polytechnique/LIX\\
and Centre national de la recherche scientifique (CNRS)\\
France
   \email{smith@lix.polytechnique.fr}}
\end{document}